\documentclass{amsart}
\usepackage{amssymb, amsmath, verbatim} 
\usepackage{enumitem}
\usepackage{multirow}
\setlist[itemize]{leftmargin=12mm}
\setlist[enumerate]{leftmargin=12mm}

\DeclareMathOperator{\ord}{\upsilon}

\newcommand{\Q}{{\mathbb Q}}
\newcommand{\Z}{{\mathbb Z}}

\newcommand{\F}{{\mathbb F}}
\newcommand{\PP}{{\mathbb P}}

\newcommand{\cA}{\mathcal{A}}
\newcommand{\cC}{\mathcal{C}}

\newcommand{\cB}{\mathcal{B}}
\newcommand{\cG}{\mathcal{G}}
\newcommand{\cE}{\mathcal{E}}
\newcommand{\cH}{\mathcal{H}}

\newcommand{\cM}{\mathcal{M}}
\newcommand{\cN}{\mathcal{N}}

\newcommand{\OO}{{\mathcal O}}
\newcommand{\ga}{\mathfrak{a}}

\newcommand{\ff}{\mathfrak{f}}

\newcommand{\fp}{\mathfrak{m}}
\newcommand{\fq}{\mathfrak{q}}
\newcommand{\mP}{\mathfrak{P}}     
\newcommand{\sS}{\mathfrak{S}}

\begin {document}

\newtheorem{thm}{Theorem}
\newtheorem{lem}{Lemma}[section]

\newtheorem{cor}[lem]{Corollary}
\newtheorem{conj}{Conjecture}

\theoremstyle{remark}

\title[]{
On the generalized Fermat equation\\ over totally real fields 
}

\author{Heline Deconinck}
\thanks{The author's research is funded by an EPSRC studentship.}

\begin{abstract}
In a recent paper, Freitas and Siksek proved an
asypmtotic version of Fermat's Last Theorem for many
totally real fields. We prove an extension of their
result to generalized Fermat equations of the form
$A x^p+B y^p+ C z^p=0$, where $A$, $B$, $C$ are
odd integers belonging to a totally real field.
\end{abstract}
\maketitle

\section{Introduction}

Let $K$ be a totally real number field and let $\OO_K$ be its ring of
integers. In \cite{FS}, Freitas and Siksek study the 
Fermat equation $a^p+b^p+c^p=0$ with $a$, $b$, $c \in \OO_K$
and $p$ prime. For now let $S$ be the set of primes of $\OO_K$
above $2$ and let $\OO_S$ be the ring of $S$-integer and $\OO_S^*$
be the group of $S$-units. Freitas and Siksek give a criterion
for the non-existence of solutions $a$, $b$, $c \in \OO_K$ with $abc \ne 0$
for $p$ sufficiently large in terms of the solutions to the 
$S$-unit equation $\lambda+\mu=1$. The proof uses modularity 
and level lowering arguments over totally real fields.
It is natural to seek an extention of the work of Freitas and Siksek to 
generalized Fermat equations $A a^p+B b^p+ C c^p=0$,
for given non-zero coefficients $A$, $B$, $C \in \OO_K$.
In this paper 
we show that the results of Freitas and Siksek can indeed be extended
to any choice of \emph{odd} coefficients 
$A$, $B$, $C$, 
provided the set $S$ is enlarged to contain
the primes dividing $ABC$ as well as the primes dividing $2$.

We now state our results precisely. 
As in \cite{FS}, our results will sometimes be conditional on the following
standard conjecture.
\begin{conj}[\lq\lq Eichler--Shimura\rq\rq]\label{conj:ES}
Let $K$ be a totally real field. Let $\ff$ be a Hilbert newform
of level $\cN$ and parallel weight $2$, and write $\Q_\ff$ for the field
generated by its eigenvalues. Suppose that $\Q_\ff=\Q$. 
Then there is an elliptic curve $E_\ff/K$ with conductor $\cN$
having the same $\mathrm{L}$-function as $\ff$.
\end{conj}
Let $A$, $B$, $C$ be non-zero elements of $\OO_K$, and let $p$ be a prime. 
Consider the equation 
\begin{equation}\label{eqn:Fermat}
Aa^p+Bb^p+Cc^p=0, \qquad a,b,c\in \OO_K; 
\end{equation}
we shall refer to this as \emph{the generalized Fermat equation over $K$ with
coefficients $A$, $B$, $C$ and exponent $p$}.
A solution $(a,b,c)$ is called \textbf{trivial}
if $abc=0$, otherwise \textbf{non-trivial}. 
The following notation shall be fixed throughout the paper.
\begin{equation}\label{eqn:ST} 
\begin{gathered}
R = \text{Rad}(ABC) =  \prod_{\substack{ \fq \mid  ABC \\ \fq \text{ prime in } K}} \fq \\
S =\{ \mP \; :\; \text{$\mP$ is a prime of $\OO_K$ such that $\mP \mid 2R$}\}\\ 
T  =\{ \mP \; :\; \text{$\mP$ is a prime of $\OO_K$ above $2$}\}, \\
U=\{ \mP \in T \; : \; f(\mP/2) = 1 \}, \qquad
V=\{ \mP \in T \; : \; 3 \nmid \ord_\mP(2) \} 
\end{gathered}
\end{equation}
Here $f(\mP/2)$ denotes the residual degree of $\mP$.
As in \cite{FS}, we need an assumption which we refer to throughout
the paper as (ES):
\[ \label{ES}
\text{\bf (ES)} \qquad
\left\{
\begin{array}{lll}
\text{either $[K:\Q]$ is odd;}\\
\text{or $U \ne \emptyset$;}\\
\text{or Conjecture~\ref{conj:ES} holds for $K$.}
\end{array}
\right.
\]
\begin{thm}\label{thm:FermatGen}
Let $K$ be a totally real  field  satisfying
 (ES). 
Let 
$A$, $B$, $C \in \OO_K$, and 
suppose that $A$, $B$, $C$ are odd,
in the sense that if $\mP \mid 2$ is a prime of $\OO_K$ then
$\mP \nmid ABC$. 
Write $\OO_{S}^*$ for the set of $S$-units of $K$.
Suppose that for every solution $(\lambda,\mu)$ to the $S$-unit equation 
\begin{equation}\label{eqn:sunit}
\lambda+\mu=1, \qquad \lambda,\, \mu \in \OO_{S}^* \, ,
\end{equation}
there is 
 \begin{enumerate}
\item[(A)] either some $\mP \in U$ that satisfies 
$\max\{ \lvert \ord_{\mP} (\lambda) \rvert, \lvert \ord_{\mP}(\mu) \rvert \} 
\le 4 \ord_{\mP}(2)$,
\item[(B)] or some $\mP \in V$ that satisfies both 
$\max\{ \lvert \ord_{\mP} (\lambda) \rvert, \lvert \ord_{\mP}(\mu) \rvert \} 
\le 4 \ord_{\mP}(2)$, and 
$\ord_{\mP}(\lambda \mu) \equiv \ord_{\mP}(2) \pmod{3}$.
\end{enumerate}
Then there is some constant $\cB=\cB(K,A,B,C)$  such that 
the generalized Fermat equation~\eqref{eqn:Fermat} with exponent $p$ and coefficients $A$, $B$, $C$
does not have non-trivial solutions with $p>\cB$.
\end{thm}
 
Theorem~\ref{thm:FermatGen} gives a bound on the exponent
of non-trivial solutions to the generalized Fermat equation \eqref{eqn:Fermat}
provided certain hypotheses are satisfied.
There are practical algorithms for determining the solutions to
$S$-unit equations (e.g.\ \cite{Smart}), so these hypotheses can always be checked
for specific $K$, $A$, $B$, $C$. The following theorem is an example
where the $S$-unit equation can still be solved, even though the coefficients
are not completely fixed. 
\begin{thm}\label{thm:d5mod8}
Let $d \geq 13$ be squarefree, satisfying $d \equiv 5 \pmod{8}$ and let $q \geq 29$ be a prime such that $q \equiv 5 \pmod{8}$ and $\left( \frac{d}{q} \right) = -1$.
Let $K=\Q(\sqrt{d})$ and assume Conjecture~\ref{conj:ES} 
for $K$.
Then there is an effectively computable constant $\cB_{K,q}$  such that
for all primes $p > \cB_{K,q}$, the Fermat equation 
$$x^p+y^p+q^rz^p=0$$ 
has no non-trivial solutions with exponent $p$.
\end{thm}

\section{Preliminaries}
We shall need the theoretical machinary of modularity, irreducibility
of Galois representations and level lowering. 
These tools and the way we use them is practically identical to
\cite{FS} which we refer the reader to for more details.

\subsection{The Frey curve and its modularity}
We shall need the following
recently proved theorem \cite{FHS}.
\begin{thm}[Freitas, Le Hung and Siksek]
\label{thm:modgen}
Let $K$ be a totally real field. Up to isomorphism over $\overline{K}$,
there are at most finitely many non-modular 
elliptic curves $E$ over $K$. Moreover, if $K$ is real quadratic,
then all elliptic curves over $K$ are modular.
\end{thm}

We shall associate to a solution $(a,b,c)$ of \eqref{eqn:Fermat}
the following Frey elliptic curve.
\begin{equation}\label{eqn:Frey}
E : Y^2 = X(X-Aa^p)(X+Bb^p) 
\end{equation}
Before applying Theorem~\ref{thm:modgen} to the Frey curve
associated to our generalized Fermat equation \eqref{eqn:Fermat}
we shall need the following lemma.
\begin{lem}\label{lem:not1}
Let $A$, $B$, $C \in \OO_K$ be odd, and suppose that
every solution $(\lambda,\mu)$
to the $S$-unit equation \eqref{eqn:sunit}
satisfies
either condition (A) or (B) of Theorem~\ref{thm:FermatGen}.
Then $(\pm 1, \pm 1, \pm 1)$ is not a solution to
equation \eqref{eqn:Fermat}.
\end{lem}
\begin{proof}
Suppose $(\pm 1, \pm 1, \pm 1)$ is a solution to \eqref{eqn:Fermat}.
By changing signs of $A$, $B$, $C$, we may suppose that $(1,1,1)$
is a solution, and therefore that
$A+B+C=0$. 
Let $\lambda=A/C$ and $\mu=B/C$.
Clearly $(\lambda,\mu)$ is a solution to the $S$-unit equation \eqref{eqn:sunit}.

Suppose first that (A) is satisfied.
Then $U \ne \emptyset$, so there is some $\mP \mid 2$ with residue
field $\F_2$. As $A$, $B$, $C$ are odd, we have $\mP \nmid ABC$.
Reducing the relation $A+B+C=0$ mod $\mP$ we obtain
$1+1+1=0$ in $\F_2$, giving a contradiction.

Suppose now that (B) holds. 
By (B) there is some $\mP \in V$
such that $\ord_\mP(\lambda \mu) \equiv \ord_\mP(2) \pmod{3}$.
However, as $A$, $B$, $C$ are odd, $\ord_\mP(\lambda \mu)=0$.
Moreover, $3 \nmid \ord_\mP(2)$ by definition of $V$. This
gives a contradiction.
\end{proof}

\begin{cor}\label{cor:Freymod}
Let $A$, $B$, $C \in \OO_K$ be odd, and suppose that
every solution $(\lambda,\mu)$
to the $S$-unit equation \eqref{eqn:sunit}
satisfies
either condition (A) or (B) of Theorem~\ref{thm:FermatGen}.
There is some (ineffective) constant $\cA=\cA(K,A,B,C)$ 
such that for any non-trivial solution 
$(a,b,c)$ of \eqref{eqn:Fermat} 
with prime exponent $p>\cA$,
the Frey curve $E$ given by \eqref{eqn:Frey} is modular.
\end{cor}
\begin{proof}
By Theorem~\ref{thm:modgen}, there are at most finitely many possible
$\overline{K}$-isomorphism classes of elliptic curves over $K$
that are non-modular. Let $j_1,\dots,j_n \in K$ be the 
$j$-invariants of these classes. 
Write $\lambda=-Bb^p/Aa^p$. The $j$-invariant of $E_{a,b,c}$ is 
\[
j(\lambda)=2^8 \cdot (\lambda^2-\lambda+1)^3 \cdot \lambda^{-2} (\lambda-1)^{-2}.
\]
Each equation $j(\lambda)=j_i$ has at most six solutions $\lambda \in K$. 
Thus there are values $\lambda_1,\dots,\lambda_m \in K$
such that if $\lambda \ne \lambda_k$ for 
all $k$ then $E$ is modular. 
If $\lambda=\lambda_k$ then 
\[
 (-b/a)^p=A\lambda_k/B, \qquad (c/a)^p=A( \lambda_k-1)/C. 
\]
This pair of equations results in a bound for $p$
unless $-b/a$ and $c/a$ are both
roots of unity. But as $K$
is real, the only roots of unity are $\pm 1$.
If $-b/a=\pm 1$ and $c/a= \pm 1$ then \eqref{eqn:Fermat}
has a solution of the form $(\pm 1,\pm 1, \pm 1)$
contradicting Lemma~\ref{lem:not1}.
This completes the proof.
\end{proof}

\subsection{Irreducibility of mod $p$ representations of elliptic curves}
To use a generalized version of level lowering, we need the 
mod $p$ Galois representation associated to the Frey elliptic 
curve to be irreducible. The following theorem of Freitas and Siksek 
\cite[Theorem 2]{FSirred}, building on earlier work of David, Momose
and Merel, is sufficient for our purpose.
\begin{thm} \label{thm:irred}
Let $K$ be a totally real field. 
There is an effective constant $\cC_K$, depending only on $K$, such that the 
following holds. 
If $p > \cC_K$ is a rational prime,
and  $E$ is an elliptic curve over $K$ which is semistable at some $\fq \mid p$,
then $\overline{\rho}_{E,p}$ is irreducible.
\end{thm}
In \cite{FSirred} the theorem is stated for Galois totally real fields $K$,
but the version stated here follows immediately on replacing $K$ by its 
Galois closure.

\subsection{Level Lowering}

As before, $K$ is a totally real field. Let $E/K$ be an elliptic curve of conductor $\cN$
and $p$ a rational prime.
For a prime ideal $\fq$ of $K$ denote by $\Delta_\fq$ the discriminant of a local
minimal model for $E$ at $\fq$.
Let
\begin{equation}\label{eqn:Np}
\cM_p := \prod_{
\substack{\fq \Vert \cN,\\ p \mid \ord_\fq(\Delta_\fq)}
} {\fq}, \qquad\quad \cN_p:=\frac{\cN}{\cM_p} \, .
\end{equation}
The ideal $\cM_p$ is precisely the product of the primes where we want to lower the level.
For a Hilbert eigenform $\ff$ over $K$, denote the field generated by its 
eigenvalues by $\Q_\ff$.
The following level-lowering recipe is derived by Freitas and 
Siksek \cite{FS} from the works of
Fujiwara, Jarvis and Rajaei.
\begin{thm}\label{thm:levell} 
With the above notation, suppose the following
\begin{enumerate}
\item[(i)] $p\ge 5$ and $p$ is unramified in $K$
\item[(ii)] $E$ is modular,
\item[(iii)] $\overline{\rho}_{E,p}$ is irreducible,
\item[(iv)] $E$ is semistable at all $\fq \mid p$,
\item[(v)]  $p \mid \ord_\fq(\Delta_\fq)$ for all $\fq \mid p$. 
\end{enumerate}
Then, there is a Hilbert eigenform $\ff$ 
of parallel weight $2$ that is new at level $\cN_p$ and some prime $\varpi$ of $\Q_\ff$
such that $\varpi \mid p$
and $\overline{\rho}_{E,p} \sim \overline{\rho}_{\ff,\varpi}$.
\end{thm}

\section{Conductor of the Frey curve}
Let $(a,b,c)$ be a non-trivial solution to the Fermat equation~\eqref{eqn:Fermat}.
Write
\begin{equation}\label{eqn:cG}
\cG_{a,b,c}=a \OO_K+b \OO_K+ c \OO_K,
\end{equation}
which we naturally think of as the greatest
common divisor of $a$, $b$, $c$. 
Over $\Q$, or over a number field of class number $1$ it
is natural to scale the solution $(a,b,c)$ so that $\cG_{a,b,c}=1 \cdot \OO_K$,
but this is not possible in general. The primes that divide all of $a$, $b$, $c$
can be additive primes for the Frey curve, and additive primes
are not removed by the 
level lowering recipe given above. To control the final level we 
need to control $\cG_{a,b,c}$. 
Following \cite{FS},
we fix a set  
\[
\cH=\{\fp_1,\dots,\fp_h\}
\]
of prime ideals $\fp_i \nmid 2 R$, which is a set of representatives
for the ideal classes of $\OO_K$.
For an non-zero ideal $\ga$ of $\OO_K$, we denote by
$[\ga]$ the class of $\ga$ in the class group. We denote
$[\cG_{a,b,c}]$ by $[a,b,c]$. The following is Lemma 3.2
of \cite{FS}, and states that we can always scale
our solution $(a,b,c)$ so that the gcd belongs to $\cH$.

\begin{lem}\label{lem:gcd}
Let $(a,b,c)$ be a non-trivial solution to \eqref{eqn:Fermat}.
There is a non-trivial integral solution $(a^\prime,b^\prime,c^\prime)$ to \eqref{eqn:Fermat} such 
that the following hold.
\begin{enumerate}
\item[(i)] For some $\xi \in K^*$,
\[
a^\prime= \xi a, \qquad b^\prime= \xi b, \qquad c^\prime=\xi c.
\]
\item[(ii)] $\cG_{a^\prime,b^\prime,c^\prime} = \fp \in \cH$.
\item[(iii)] $[a^\prime,b^\prime,c^\prime]=[a,b,c]$.
\end{enumerate}
\end{lem}

\begin{lem}\label{lem:cond}

Let $(a,b,c)$ be a non-trivial solution to the 
Fermat equation \eqref{eqn:Fermat} with odd prime exponent $p$,
and scaled as in Lemma~\ref{lem:gcd} so that
$\cG_{a,b,c}=\fp \in \cH$. 

Write $E=E_{a,b,c}$
for the Frey curve in \eqref{eqn:Frey}, and let $\Delta$ be its discriminant. 
For a prime $\fq$ we write $\Delta_\fq$ for the minimal discriminant at $\fq$.
Then at all $\fq \notin S \cup \{\fp\}$, 
the model $E$ is minimal, semistable, and satisfies 
$p \mid \ord_\fq(\Delta_\fq)$.
Let $\cN$ be the conductor of $E$, and let $\cN_p$
be as defined in \eqref{eqn:Np}. Then
\begin{equation}\label{eqn:cnp}
\cN=\fp^{s_{\fp}} 
\cdot
\prod_{\mP \in S} \mP^{r_\mP}
\cdot 
\prod_{\substack{\fq \mid abc \\ 
\fq \notin S \cup \{\fp\}
}} \fq \, , \qquad \qquad
\cN_p= \fp^{s^\prime_{\fp}}  \cdot
\prod_{\mP \in S} \mP^{r^\prime_\mP},
\end{equation}
where $0 \le r^\prime_\mP \le r_\mP \le 
2 + 6\ord_{\mP} (2)$ for  $\mP \mid 2$,
and 
$0 \le r^\prime_\mP \le r_\mP \le 
2$ for  $\mP \mid R$, and
$0 \leq s^\prime_{\fp} \leq s_{\fp} \leq 2$.
\end{lem}

\begin{proof}
The discriminant of the model given by $E$ is $16 (ABC)^2 (abc)^{2p}$,
thus the primes appearing in $\cN$ will be either primes dividing $2R$
or dividing $abc$. 
For $\mP \mid 2$ we have
$r_\mP=\ord_\mP(\cN) \le 2+6 \ord_{\mP}(2)$ by \cite[Theorem IV.10.4]{SilvermanII}; this proves the correctness of the bounds for the
exponents in $\cN$ and $\cN_p$ at even primes, and we will restrict
our attention to odd primes.
As $E$ has
full 2-torsion over $K$, the wild part of the conductor of $E/K$ vanishes 
(\cite{SilvermanII}, page 380) at all odd $\fq$, and so $\ord_\fq(\cN_p)
\le \ord_\fq(\cN) \le 2$. This proves the correctness of the bounds
for the exponents in $\cN$ and $\cN_p$ at $\fq$ that divide $R$
and for $\fq=\fp$. 

It remains to consider $\fq \mid abc$
satisfying $\fq \not \in S \cup \{ \fp\}$. It is easily checked
that the model \eqref{eqn:Frey} is minimal and has multiplicative
reduction at such $\fq$, and it is therefore clear that
$p \mid \ord_\fq(\Delta)=\ord_\fq(\Delta_\fq)$. It follows
that $\ord_\fq(\cN)=1$ and, from the recipe for 
$\cN_p$ in \eqref{eqn:Np} that $\ord_\fq(\cN_p)=0$.
This completes the proof.

\end{proof}

\section{Level Lowering for the Frey Curve}

\begin{thm}\label{thm:ll2}
Let $K$ be a totally real field  satisfying (ES).
Let $A$, $B$, $C \in \OO_K$ be odd, and suppose that
every solution $(\lambda,\mu)$
to the $S$-unit equation \eqref{eqn:sunit}
satisfies
either condition (A) or (B) of Theorem~\ref{thm:FermatGen}.
There is a constant 
$\cB=\cB(K,A,B,C)$ 
depending only on $K$ and $A$, $B$, $C$ such that 
the following hold.
Let $(a,b,c)$
be a non-trivial solution to the 
generalized Fermat equation \eqref{eqn:Fermat} 
with prime exponent $p>\cB$,
and rescale $(a,b,c)$ as in Lemma~\ref{lem:gcd} so that it remains integral and
satisfies $\cG_{a,b,c}=\fp$ for some $\fp \in \cH$.
Write $E=E_{a,b,c}$ for the Frey curve given in \eqref{eqn:Frey}.
Then there is an elliptic curve $E^\prime$ over $K$ such that
\begin{enumerate}
\item[(i)] the conductor of $E^\prime$ is divisible only by primes in $S \cup \{ \fp \}$;
\item[(ii)] $\# E^\prime(K)[2]=4$;
\item[(iii)] $\overline{\rho}_{E,p} \sim \overline{\rho}_{E^\prime,p}$;
\end{enumerate}
Write $j^\prime$ for the $j$-invariant of $E^\prime$.
Then,
\begin{enumerate}
\item[(a)] for $\mP \in U$, we have $\ord_\mP(j^\prime)<0$;
\item[(b)] for $\mP \in V$, we have either $\ord_\mP(j^\prime)<0$
or $3 \nmid \ord_\mP(j^\prime)$;
\item[(c)] for $\fq \notin S$, we have $\ord_\fq(j^\prime) \ge 0$.
\end{enumerate}
In particular, $E^\prime$ has potentially good reduction away from $S$. 
\end{thm}

\begin{proof}
We first observe, by Lemma~\ref{lem:cond}, that
$E$ is semistable outside $S \cup \{\fp\}$. By taking $\cB$ 
to be sufficiently large, we see
from Corollary~\ref{cor:Freymod} that $E$ is modular,
and from Theorem~\ref{thm:irred} that $\overline{\rho}_{E,p}$
is irreducible. 
Applying Theorem~\ref{thm:levell}
and Lemma~\ref{lem:cond} we see that $\overline{\rho}_{E,p} \sim \overline{\rho}_{\ff,\varpi}$
for a Hilbert newform $\ff$ of level $\cN_p$
and some prime $\varpi \mid p$ of $\Q_\ff$. 
Here $\Q_\ff$ is the field
generated by the Hecke eigenvalues of $\ff$. 
The remainder of the proof is identical to the proof
of \cite[Theorem 9]{FS}, and so we omit
 the details, except that we point out that it is here
that we make use of assumption (ES).

\end{proof}

The constant $\cB$ is ineffective as it depends on the ineffective
constant $\cA$ in Corollary~\ref{cor:Freymod}. However, if $K$
is a real quadratic field then we do not need that corollary
as we know modularity from Theorem~\ref{thm:modgen}. In this 
case the arguments of \cite{FS} produce an effective constant $\cB$.

\section{Elliptic curves with full $2$-torsion and solutions to the $S$-unit equation}
Theorem~\ref{thm:ll2} relates non-trivial solutions
of the Fermat equation to elliptic curves with full $2$-torsion
having potentially good reduction outside $S$. 
There is a well-known correspondence between 
such elliptic curves and solutions of 
the $S$-unit equation \eqref{eqn:sunit}
that we now sketch. 

Consider an elliptic curve over $K$ with full $2$-torsion,
\begin{equation}\label{eqn:e123}
y^2=(x-a_1)(x-a_2)(x-a_3).
\end{equation}
where $a_1$, $a_2$, $a_3$ are distinct. The
\textbf{cross ratio}
\[
\lambda=\frac{a_3-a_1}{a_2-a_1}
\]
belongs to $\PP^1(K)-\{0,1,\infty\}$.
Moreover, any $\lambda \in \PP^1(K)-\{0,1,\infty\}$
can be written as a cross ratio of three distinct
$a_1$, $a_2$, $a_3$ in $K$ and hence comes from
an elliptic curve with full $2$-torsion.
Write $\sS_3$ for the symmetric group on $3$ letters.
The action of $\sS_3$
on the triple $(e_1, e_2, e_3)$ extends via the cross ratio
in a well-defined manner to an action on $\PP^1(K)-\{0,1,\infty\}$.
The orbit of $\lambda \in \PP^1(K)-\{0,1,\infty\}$ under the action of
$\sS_3$
is
\begin{equation}\label{eqn:orbit}
\left\{ \lambda, 
\frac{1}{\lambda},
1-\lambda, 
\frac{1}{1-\lambda}, 
\frac{\lambda}{\lambda-1}, 
\frac{\lambda-1}{\lambda}
\right\}.
\end{equation}
It follows from the theory of Legendre elliptic curves
(\cite[Pages 53--55]{SilvermanI})
that the cross ratio in fact defines 
a bijection between elliptic curves over $K$
having full $2$-torsion (up to isomorphism over $\overline{K}$),
and $\lambda$-invariants up to the action of $\sS_3$. 
Under this bijection,
the $\sS_3$-orbit of a given $\lambda \in \PP^1(K)\backslash\{0,1,\infty\}$
is associated to the $\overline{K}$-isomorphism
class of the \textbf{Legendre elliptic curve}
$y^2=x(x-1)(x-\lambda)$. We would like to understand the $\lambda$-invariants
that correspond to elliptic curves over $K$ with full $2$-torsion
and potentially good reduction outside $S$. The $j$-invariant
of the Legendre elliptic curve is given by
\begin{equation}\label{eqn:j}
j(\lambda)=2^8 \cdot \frac{(\lambda^2-\lambda+1)^3}{\lambda^2 (1-\lambda)^2} \, .
\end{equation}
The Legendre elliptic curve (and therefore its $\overline{K}$-isomorphism
class) has potentially good reduction outside $S$ if and only if 
$j(\lambda)$ belongs to $\OO_S$. It easily follows from \eqref{eqn:j}
that this happens precisely when both $\lambda$ and $1-\lambda$
are $S$-units (recall that $S$ includes all the primes above $2$); 
in other words, this is equivalent to $(\lambda,\mu)$
being a solution to the $S$-unit equation \eqref{eqn:sunit},
where $\mu=1-\lambda$. 
Let
$\Lambda_{S}$ be the set of solutions to the $S$-unit
equation \eqref{eqn:sunit}:
\begin{equation}\label{eqn:LambdaS}
\Lambda_{S}=\{(\lambda,\mu)\; : \; \lambda+\mu=1, \qquad \lambda,\; \mu \in \OO_{S}^*\}.
\end{equation}
It is easy to see that 
the action of $\sS_3$ on $\PP^1(K)-\{0,1,\infty\}$ induces
a well-defined action on $\Lambda_{S}$
given by
\[
(\lambda,\mu)^{\sigma} =(\lambda^{\sigma},1-\lambda^\sigma).
\]
We denote by $\sS_3 \backslash \Lambda_{S}$ the set
of $\sS_3$-orbits in $\Lambda_{S}$. We deduce the following.
\begin{lem}\label{lem:prebij}
Let
$\cE_{S}$ be set of all elliptic curves over $K$ with 
full $2$-torsion and
potentially good
reduction outside $S$.
 Define the equivalence relation $E_1 \sim E_2$
on $\cE_{S}$ to mean that $E_1$ and $E_2$ are isomorphic
over $\overline{K}$. 
There is a well-defined bijection
\[
\Phi \; : \; \cE_{S}/\sim \; \longrightarrow \; \sS_3 \backslash \Lambda_{S}
\]
which sends the class of an elliptic curve
given by \eqref{eqn:e123} to the orbit of
\[
\left(\frac{a_3-a_1}{a_2-a_1},\; \frac{a_2-a_3}{a_2-a_1} \right)
\]
in $\sS_3 \backslash \Lambda_{S}$; the map $\Phi^{-1}$ sends 
the orbit of $(\lambda,\mu)$ to the class of the Legendre
elliptic curve $y^2=x(x-1)(x-\lambda)$. 
\end{lem}

We shall need the following for the proof of Theorem~\ref{thm:FermatGen}.
\begin{lem}\label{lem:jval}
Let $E^\prime \in \cE_{S}$ and 
suppose that its $\sim$-equivalence class 
corresponds via $\Phi$ to the orbit of $(\lambda,\mu) \in \Lambda_{S}$.
Let $j^\prime$ be the $j$-invariant of $E^\prime$ and $\mP \in T$. 
Then
\begin{itemize}
\item[(i)] $\ord_\mP(j^\prime) \ge 0$ if and only if
$\max\{\lvert \ord_\mP(\lambda) \rvert,
\lvert \ord_\mP(\mu) \rvert \} 
 \le 4 \ord_\mP(2)$,
\item[(ii)] $3 \mid \ord_\mP(j^\prime)$ if and only 
$\ord_\mP(\lambda \mu) \equiv \ord_\mP(2) \pmod{3}$.
\end{itemize}
\end{lem}
\begin{proof}
Observe that
\begin{equation}\label{eqn:jlambda}
j^\prime=
j(\lambda)=2^8 \cdot \frac{(\lambda^2-\lambda+1)^3}{\lambda^2 (\lambda-1)^2}
=2^8 \cdot \frac{(1-\lambda \mu)^3}{(\lambda \mu)^2} \, .
\end{equation}
From this we immediately deduce (ii). 
Let 
\[
m=\ord_\mP(\lambda), \qquad n=\ord_\mP(\mu), \qquad
t=\max(\lvert m \rvert, \lvert n \rvert).
\]
If $t=0$ then $\ord_\mP(j^\prime) \ge 8 \ord_\mP(2) > 0$,
and so (i) holds. We may therefore suppose that $t>0$.
Now the relation $\lambda+\mu=1$
forces either $m=n=-t$, or $m=0$ and $n=t$, or $m=t$ and $n=0$.
Thus $\ord_\mP(\lambda \mu)=-2t<0$ or $\ord_\mP(\lambda \mu)=t>0$.
In either case, from~\eqref{eqn:j},
\[
\ord_\mP(j^\prime)=8 \ord_\mP(2)-2 t.
\]
This proves (i).
\end{proof}

\section{Proof of Theorem~\ref{thm:FermatGen}}
Let $K$ be a totally real field satisfying assumption (ES).
Let $S$, $T$, $U$, $V$ be as in \eqref{eqn:ST}.
Let $\cB$ be as in Theorem~\ref{thm:ll2}, and let $(a,b,c)$
be a non-trivial solution to the Fermat equation~\eqref{eqn:Fermat}
with exponent $p>\cB$,
scaled so that 
$\cG_{a,b,c}=\fp$ with $\fp \in \cH$.
Applying
Theorem~\ref{thm:ll2} gives an elliptic curve $E^\prime/K$
with full $2$-torsion and potentially good reduction outside $S$
whose $j$-invariant $j^\prime$ satisfies:
\begin{enumerate}
\item[(a)] for all $\mP \in U$, we have $\ord_\mP(j^\prime) <0$;
\item[(b)] for all $\mP \in V$, we have $\ord_\mP(j^\prime)<0$
or $3 \nmid \ord_\mP(j^\prime)$.
\end{enumerate}
Let $(\lambda,\mu)$ be a solution to $S$-unit equation \eqref{eqn:sunit},
whose $\sS_3$-orbit corresponds to the $\overline{K}$-isomorphism
class of $E^\prime$ as in Lemma~\ref{lem:prebij}.  
By Lemma~\ref{lem:jval} and (a), (b) we know that
\begin{enumerate}
\item[($a^\prime$)] for all $\mP \in U$, we have 
$\max\{
\lvert \ord_\mP(\lambda)
\rvert, \lvert \ord_\mP(\mu) \rvert \} > 4 \ord_\mP(2)$;
\item[($b^\prime$)] for all $\mP \in V$, we have
$\max\{
\lvert \ord_\mP(\lambda)
\rvert, \lvert \ord_\mP(\mu) \rvert \} > 4 \ord_\mP(2)$
or $\ord_\mP(\lambda \mu) \not\equiv \ord_\mP(2) \pmod{3}$.
\end{enumerate}
These contradict assumptions (A) and (B) of Theorem~\ref{thm:FermatGen},
completing the proof.

\section{The $S$-unit equation over real quadratic fields}
To prove Theorem ~\ref{thm:d5mod8}
we need to understand the solutions to
the $S$-unit equation \eqref{eqn:sunit} for
real quadratic fields $K$. This is easier when
$S$ is small in size.
\begin{lem}\label{lem:makeintegral}
Suppose $\lvert S \rvert=2$. Let $(\lambda,\mu) \in \Lambda_S$. 
Then, there is some element $\sigma \in \sS_3$
so that $(\lambda^\prime,\mu^\prime)=(\lambda,\mu)^\sigma$
satisfies $\lambda^\prime$, $\mu^\prime \in \OO_K$.
\end{lem}
\begin{proof}
As $\mu^\prime=1-\lambda^\prime$ we need only
find some element $\sigma \in \sS_3$ so that
 $\lambda^\prime=\lambda^\sigma \in \OO_K$.
Write $S=\{\mP_1,\mP_2\}$. 
If $\ord_{\mP_i}(\lambda) \ne 0$ for $i=1$, $2$, 
then let $\lambda^\prime=\lambda/(\lambda-1)$,
which will have non-negative valuation at $\mP_i$
and so belongs to $\OO_K$. Thus without loss of
generality we may suppose that $\ord_{\mP_1}(\lambda)=0$.
Now if $\ord_{\mP_2}(\lambda) \ge 0$ then $\lambda^\prime=\lambda \in \OO_K$,
and if $\ord_{\mP_2}(\lambda)<0$ then $\lambda^\prime=1/\lambda \in \OO_K$. 
\end{proof}

For the remainder of this section $d$ denotes a squarefree integer $\geq 13$
that satisfies $d \equiv 5 \pmod{8}$ and $q\geq 29$ a prime satisfying $q \equiv 5 \pmod{8}$ and
$\left( \frac{d}{q} \right)=-1$.
Let $K$ denotes the real quadratic field $\Q(\sqrt{d})$.
It follows that both $q$ and $2$
are inert in $K$. We let $S=\{2,q\}$.
\begin{lem}\label{lem:makeint}
Let $K$ and $S$ be as above, and let $(\lambda,\mu) \in \Lambda_S$.
Then $\lambda$, $\mu \in \Q$ if and only if $(\lambda,\mu)$
belongs to the $\sS_3$-orbit $\{(1/2,1/2),\; (2,-1),\; (-1,2)\}
\subseteq \Lambda_S$.
\end{lem}
\begin{proof}
Suppose $\lambda$, $\mu \in \Q$. By Lemma~\ref{lem:makeintegral}
we may suppose that $\lambda$ and $\mu$ belong to 
$\OO_K \cap \Q=\Z$ and hence 
$\lambda=\pm 2^{r_1} q^{s_1}$, $\mu=\pm 2^{r_2} q^{s_2}$
where $r_i \ge 0$ and $s_i \ge 0$. As $\lambda+\mu=1$
we see that one of $r_1$, $r_2$ is $0$ and likewise one of 
$s_1$, $s_2=0$. Without loss of generality $r_2=0$. If $s_2=0$
too then we have $\lambda\pm 1=1$ which forces $(\lambda,\mu)=(2,-1)$
as required. We may therefore suppose that $s_1=0$.
Hence $\pm 2^{r_1} \pm q^{s_2}=1$. If $s_2=0$ then again
we obtain $(\lambda,\mu)=(2,-1)$, so suppose $s_2>0$. 

We now easily check that $r_1=1$ and $r_1=2$ are both incompatible
with our hypotheses on $q$. Thus $r_1 \ge 3$ and so $\mu=\pm q^{s_2} \equiv 1 \pmod{8}$.
As $q \equiv 5 \pmod{8}$, we have $\mu=q^{2t}$
for some integer $t \ge 1$. Hence $(q^t+1)(q^t-1)=\mu-1=-\lambda=\mp 2^{r_1}$.
This implies that $q^t+1=2^{a}$ and $q^t-1=2^{b}$ where $a \ge b \ge 1$. Subtracting we
have $2^a-2^b=2$ and so $b=1$ and $q=3$ giving a contradiction.
\end{proof}

We follow \cite{FS} in calling the elements of the orbit
$\{(1/2,1/2),\; (2,-1),\; (-1,2)\}$ \textbf{irrelevant},
and in calling other elements of $\Lambda_S$ relevant.
Next we give a parametrization of all relevant elements 
of $\Lambda_S$. This the analogue of \cite[Lemma 6.4]{FS},
and shows that such a parametrization is possible even
though our set $S$ is larger, containing the odd prime $q$.

\begin{lem}\label{lem:param}
Up to the action of $\sS_3$, every relevant $(\lambda,\mu) \in \Lambda_{S}$
has the form 
\begin{equation}\label{eqn:parasol}
\lambda=\frac{\eta_1 \cdot 2^{2r_1} \cdot q^{2s_1} -\eta_2 \cdot  q^{2s_2}+1 +v \sqrt{d}}{2},
\qquad 
\mu=\frac{\eta_2 \cdot q^{2s_2} -\eta_1 \cdot 2^{2r_1} \cdot q^{2s_1}+1 -v \sqrt{d}}{2}
\end{equation}
where
\begin{equation}\label{eqn:paracond1}
\eta_1=\pm 1, \qquad \eta_2=\pm 1, \qquad r_1  \ge 0, \qquad s_1, s_2 \geq 0, \qquad s_1 \cdot s_2=0, \qquad v \in \Z, \qquad v \ne 0
\end{equation}
are related by
\begin{gather}\label{eqn:paracond2}
(\eta_1 \cdot 2^{2 r_1} \cdot q^{2s_1}-\eta_2 \cdot q^{2s_2}+1)^2 -d v^2   =\eta_1 \cdot 2^{2r_1+2} \cdot q^{2s_1}\\
\label{eqn:paracond3}
(\eta_2 \cdot q^{2s_2}-\eta_1 \cdot 2^{2r_1} \cdot q^{2s_1}+1)^2 -d v^2   =\eta_2 \cdot 2^{2} \cdot q^{2s_2}.
\end{gather}
\end{lem}
\begin{proof}
If $\eta_1$, $\eta_2$, $r_1$, $s_1$, $s_2$ and $v$ satisfy \eqref{eqn:paracond1}, \eqref{eqn:paracond2},
\eqref{eqn:paracond3}
and $\lambda$, $\mu$ are given by \eqref{eqn:parasol}, 
it is clear that $(\lambda,\mu)$
is a relevant element of $\Lambda_S$.

Conversely, suppose $(\lambda,\mu)$ is a relevant element of $\Lambda_S$. 
By Lemma~\ref{lem:makeint}, we may suppose that $\lambda$, $\mu \in \OO_K$,
and that $\lambda$, $\mu \notin \Q$.
As $S=\{2,q\}$ we can write $\lambda=2^{r_1} q^{s_1} \lambda^\prime$ and $\mu=2^{r_2} q^{s_2} \mu^\prime$
where $\lambda^\prime$ and $\mu^\prime$ are units.  As $\lambda+\mu=1$ we have
$r_1 r_2=0$ and $s_1 s_2=0$. Swapping $\lambda$ and $\mu$ if necessary, we can suppose that 
$r_2=0$.
Let $x \mapsto \overline{x}$ denote conjugation in $K$.
Then
\[
\lambda \overline{\lambda}=\eta_1 \cdot 2^{2 r_1} \cdot q^{2s_1}, \qquad \mu \overline{\mu}=\eta_2 \cdot q^{2s_2}, \qquad
\eta_1=\pm 1, \qquad \eta_2=\pm 1.
\]

Now, 
\begin{equation*}
\lambda+\overline{\lambda} = \lambda \overline{\lambda} - (1-\lambda)(1-\overline{\lambda}) +1
						 = \lambda \overline{\lambda} - \mu \overline{\mu} +1 
						 = \eta_1 \cdot 2^{2r_1} \cdot q^{2s_1} - \eta_2 \cdot q^{2s_2} +1 \, .
\end{equation*}
Moreover we can write $\lambda-\overline{\lambda}=v \sqrt{d}$, where $v \in \Z$, and as $\lambda \notin \Q$,
we have $v \ne 0$. The expressions for $\lambda+\overline{\lambda}$ and $\lambda-\overline{\lambda}$
give the expression for $\lambda$ in \eqref{eqn:parasol}, and we deduce the expression for $\mu$
from $\mu=1-\lambda$. Finally, \eqref{eqn:paracond2} follows from the identity
\[
(\lambda+\overline{\lambda})^2-(\lambda-\overline{\lambda})^2=4 \lambda \overline{\lambda},
\]
and \eqref{eqn:paracond3} from the corresponding identity for $\mu$.
\end{proof}

\begin{lem}\label{lem:nottable}
Let $d  \equiv 5 \pmod{8}$ be squarefree $d \geq 13$ and $q \geq 29$ a prime
such that $q \equiv 5 \pmod{8}$ and $\left( \frac{d}{q} \right) =  -1$.  Then
there are no relevant elements of $\Lambda_{S}$. 
\end{lem}

\begin{proof}
We apply Lemma~\ref{lem:param}. In particular, $s_1 s_2=0$.
Suppose first that $s_1>0$. 
Thus $s_2=0$.
As $( d/q )=-1$,
we have from \eqref{eqn:paracond2} that $q^{s_1} \mid v$ and 
$q^{s_1} \mid  (\eta_2 -1)$. Hence $\eta_2=1$. 
Now \eqref{eqn:paracond2} can be rewritten
as
\[
2^{4 r_1} q^{2s_1} - d (v/q^{s_1})^2=\eta_1 2^{2 r_1+2} \, .
\]
Thus $(d/q)=(-\eta_1/q)=1$ as $q \equiv 5 \pmod{8}$.
This is a contradiction.

\medskip

Thus, henceforth, $s_1=0$.
Next suppose that $s_2=0$. We will consider the subcases $\eta_2=-1$
and $\eta_2=1$ separately and obtain contradictions in both 
subcases showing that $s_2>0$. 
Suppose $\eta_2 = -1$. From \eqref{eqn:paracond3} we have 
$2^{4r_1}  - dv^2 = - 4$. If $r_1=0$ or $1$
then $d=5$ and if $r_1 \ge 2$ then $d \equiv 1 \pmod{8}$,
giving a contradiction.
Hence suppose $\eta_2=1$.
From  \eqref{eqn:paracond2}, we have
$2^{4r_1}-dv^2 = \eta_1 2^{2r_1+2}$.
If $r_1=0$, $1$, $2$ then $dv^2=1 \pm 4$, $dv^2= 16 \pm 16$,
$dv^2= 256 \pm 64$ all of which contradict the assumptions
on $d$ or the fact that $v \ne 0$ (by \eqref{eqn:paracond1}).
If $r_1 \ge 3$ then $2^{2r_1-2}-\eta_1=d(v/2^{r_1+1})^2$
which forces $d \equiv \pm 1 \pmod{8}$, a contradiction.

\medskip

We are reduced to $s_1=0$ and $s_2>0$. From \eqref{eqn:paracond3},
as $(d/q)=-1$, we have $q^{s_2} \mid v$ and
\begin{equation}\label{eqn:divide}
q^{s_2} \mid (\eta_1 2^{2r_1}-1).
\end{equation}
The conditions $q \ge 29$ and $q \equiv 5 \pmod{8}$ force
$r_1 \ge 5$. Write $v=2^t w$ where $2 \nmid w$.
Suppose $t \le r_1-1$. From \eqref{eqn:paracond2} we have
$\eta_1 2^{2 r_1}-\eta_2 q^{2s_2}+1=2^t w^\prime$
where $2 \nmid w^\prime$. Thus ${w^\prime}^2-d w^2 \equiv 0 \pmod{8}$,
contradicting $d \equiv 5 \pmod{8}$. We may therefore suppose $t \ge r_1$.
Hence $2^{r_1} \mid (\eta_2 q^{2s_2}-1)$.
Thus $\eta_2=1$. Therefore $2^{r_1} \mid (q^{s_2}-1)(q^{s_2}+1)$.
Since $q \equiv 5 \pmod{8}$, we have $2 \mid\mid (q^{s_2}+1)$
and so
\[
2^{r_1-1} \mid (q^{s_2}-1) .
\]
As $q \equiv 5 \pmod{8}$ and $r_1 \ge 5$, we see that $s_2$ must be even,
and that $2^{r_1-2} \mid (q^{s_2/2} -1)$.
We can write $q^{s_2/2}=k \cdot 2^{r_1-2}+1$.
From \eqref{eqn:divide},
\[
k^2 2^{2r_1-4}+k 2^{r_1-1}+1=q^{s_2} \le 2^{2 r_1}+1.
\]
Hence $k=1$, $2$ or $3$. Moreover, as $q^{s_2/2} \equiv 1 \pmod{8}$,
we have $4 \mid s_2$. Hence
\[
(q^{s_2/4}-1)(q^{s_2/4}+1)=k 2^{r_1-2}.
\]
Again as $q \equiv 5 \pmod{8}$ we have $2 \mid \mid (q^{s_2/4}+1)$
and so $q^{s_2/4}+1=2$ or $6$, both of which are impossible.
This completes the proof.
\end{proof}

\section{Proof of Theorem~\ref{thm:d5mod8}}
We apply Theorem~\ref{thm:FermatGen}. By Lemma~\ref{lem:nottable}
all solutions to \eqref{eqn:sunit} are irrelevant, and the irrelevant solutions
satisfy condition (A) of Theorem~\ref{thm:FermatGen}. This
completes the proof of Theorem~\ref{thm:d5mod8}.

\end{document}